\newif\ifpagetitre            \pagetitretrue
\newtoks\hautpagetitre        \hautpagetitre={\hfil}
\newtoks\baspagetitre         \baspagetitre={\hfil}
\newtoks\auteurcourant        \auteurcourant={\hfil}
\newtoks\titrecourant         \titrecourant={\hfil}

\newtoks\hautpagegauche       \newtoks\hautpagedroite
\hautpagegauche={\hfil\the\auteurcourant\hfil}
\hautpagedroite={\hfil\the\titrecourant\hfil}

\newtoks\baspagegauche \baspagegauche={\hfil\tenrm\folio\hfil}
\newtoks\baspagedroite \baspagedroite={\hfil\tenrm\folio\hfil}

\headline={\ifpagetitre\the\hautpagetitre
\else\ifodd\pageno\the\hautpagedroite
\else\the\hautpagegauche\fi\fi}

\footline={\ifpagetitre\the\baspagetitre
\global\pagetitrefalse
\else\ifodd\pageno\the\baspagedroite
\else\the\baspagegauche\fi\fi}

\vsize=9.0in\voffset=1cm
\looseness=2


\message{fonts,}

\font\tenrm=cmr10
\font\ninerm=cmr9
\font\eightrm=cmr8
\font\teni=cmmi10
\font\ninei=cmmi9
\font\eighti=cmmi8
\font\ninesy=cmsy9
\font\tensy=cmsy10
\font\eightsy=cmsy8
\font\tenbf=cmbx10
\font\ninebf=cmbx9
\font\tentt=cmtt10
\font\ninett=cmtt9

\font\ninesl=cmsl9
\font\eightsl=cmsl8

\font\nineit=cmti9
\font\eightit=cmti8

\skewchar\ninei='177 \skewchar\eighti='177
\skewchar\ninesy='60 \skewchar\eightsy='60

\def\eightpoint{\def\rm{\fam0\eightrm} 
\normalbaselineskip=9pt
\normallineskiplimit=-1pt
\normallineskip=0pt

\textfont0=\eightrm \scriptfont0=\sevenrm \scriptscriptfont0=\fiverm
\textfont1=\ninei \scriptfont1=\seveni \scriptscriptfont1=\fivei
\textfont2=\ninesy \scriptfont2=\sevensy \scriptscriptfont2=\fivesy
\textfont3=\tenex \scriptfont3=\tenex \scriptscriptfont3=\tenex
\textfont\itfam=\eightit  \def\it{\fam\itfam\eightit} 
\textfont\slfam=\eightsl \def\sl{\fam\slfam\eightsl} 

\setbox\strutbox=\hbox{\vrule height6pt depth2pt width0pt}%
\normalbaselines \rm}

\def\ninepoint{\def\rm{\fam0\ninerm} 
\textfont0=\ninerm \scriptfont0=\sevenrm \scriptscriptfont0=\fiverm
\textfont1=\ninei \scriptfont1=\seveni \scriptscriptfont1=\fivei
\textfont2=\ninesy \scriptfont2=\sevensy \scriptscriptfont2=\fivesy
\textfont3=\tenex \scriptfont3=\tenex \scriptscriptfont3=\tenex
\textfont\itfam=\nineit  \def\it{\fam\itfam\nineit} 
\textfont\slfam=\ninesl \def\sl{\fam\slfam\ninesl} 
\textfont\bffam=\ninebf \scriptfont\bffam=\sevenbf
\scriptscriptfont\bffam=\fivebf \def\bf{\fam\bffam\ninebf} 
\textfont\ttfam=\ninett \def\tt{\fam\ttfam\ninett} 

\normalbaselineskip=11pt
\setbox\strutbox=\hbox{\vrule height8pt depth3pt width0pt}%
\let \smc=\sevenrm \let\big=\ninebig \normalbaselines
\parindent=1em
\rm}

\def\tenpoint{\def\rm{\fam0\tenrm} 
\textfont0=\tenrm \scriptfont0=\ninerm \scriptscriptfont0=\fiverm
\textfont1=\teni \scriptfont1=\seveni \scriptscriptfont1=\fivei
\textfont2=\tensy \scriptfont2=\sevensy \scriptscriptfont2=\fivesy
\textfont3=\tenex \scriptfont3=\tenex \scriptscriptfont3=\tenex
\textfont\itfam=\nineit  \def\it{\fam\itfam\nineit} 
\textfont\slfam=\ninesl \def\sl{\fam\slfam\ninesl} 
\textfont\bffam=\ninebf \scriptfont\bffam=\sevenbf
\scriptscriptfont\bffam=\fivebf \def\bf{\fam\bffam\tenbf} 
\textfont\ttfam=\tentt \def\tt{\fam\ttfam\tentt} 

\normalbaselineskip=11pt
\setbox\strutbox=\hbox{\vrule height8pt depth3pt width0pt}%
\let \smc=\sevenrm \let\big=\ninebig \normalbaselines
\parindent=1em
\rm}

\message{fin format jgr}

\hautpagegauche={\hfill\ninerm\the\auteurcourant}
\hautpagedroite={\ninerm\the\titrecourant\hfill}
\auteurcourant={R.G.\ Novikov}
\titrecourant={New global stability estimates for the Gel'fand-Calderon
inverse problem}

\magnification=1200
\font\Bbb=msbm10
\def\R{\hbox{\Bbb R}}
\def\C{\hbox{\Bbb C}}
\def\N{\hbox{\Bbb N}}
\def\S{\hbox{\Bbb S}}
\def\pa{\partial}
\def\b{\backslash}
\def\v{\varphi}

\vskip 2 mm
\centerline{\bf New global stability estimates for the Gel'fand-Calderon
inverse problem}
\vskip 2 mm
\centerline{\bf R.G.\ Novikov}
\vskip 2 mm

\noindent
{\ninerm CNRS (UMR 7641), Centre de Math\'ematiques Appliqu\'ees, Ecole
Polytechnique,}

\noindent
{\ninerm 91128 Palaiseau, France }

\noindent
{\ninerm e-mail: novikov@cmap.polytechnique.fr}

\vskip 2 mm
{\bf Abstract.}
We prove new global stability estimates for the Gel'fand-Calderon inverse
problem in 3D.
For sufficiently regular potentials this result of the present work is a
principal improvement of the result of [G.Alessandrini, Stable determination
of conductivity  by boundary measurements, Appl. Anal. {\bf 27} (1988),
153-172].

\vskip 2 mm
{\bf 1. Introduction}

We consider the equation
$$-\Delta\psi+v(x)\psi=0,\ \ x\in D,\eqno(1.1)$$
where
$$D\ \ {\rm is\ an\ open\ bounded\ domain\ in}\ \  \R^d,\ \ d\ge 2,\
\pa D\in C^2,\ v\in L^{\infty}(D).\eqno(1.2)$$
Equation (1.1) arises, in particular, in quantum mechanics, acoustics,
electrodynamics. Formally, (1.1) looks as the Schr\"odinger equation with
potential $v$ at zero energy.

We consider the map $\Phi$ such that
$${\pa\psi\over \pa\nu}\big|_{\pa D}=\Phi(\psi\big|_{\pa D}) \eqno(1.3)$$
for all sufficiently regular solutions $\psi$ of (1.1) in
$\bar D=D\cup\pa D$, where $\nu$ is the outward normal to $\pa D$. Here we
assume also that
$$0\ \ {\rm is\ not\ a\ Dirichlet\ eigenvalue\ for\ the\ operator}\ \
-\Delta+v\ \ {\rm in}\ \ D.\eqno(1.4)$$
The map $\Phi$ is called the Dirichlet-to-Neumann map for equation (1.1) and
is considered as boundary measurements for (physical model described by)
(1.1).

We consider the following inverse boundary value problem for equation (1.1):

{\bf Problem 1.1.} Given $\Phi$, find $v$.

This problem can be considered as the Gel'fand inverse boundary value problem
for the Schr\"odinger equation at zero energy (see [9], [16]). This problem
can be also considered as a generalization of the Calderon problem of the
electrical impedance tomography (see [5], [23], [16]).

We recall that the simplest interpretation of $D$, $v$ and $\Phi$ in the
framework of the electrical impedance tomography consists in the following
(see, for example, [16], [14], [13]):
$D$ is a body with isotropic conductivity $\sigma(x)$ (where
$\sigma\ge\sigma_{min}>0$),
$$\eqalignno{
&v(x)=(\sigma(x))^{-1/2}\Delta(\sigma(x))^{1/2},\ \ x\in D,&(1.5)\cr
&\Phi=\sigma^{-1/2}\bigl(\Lambda\sigma^{-1/2}+{\pa\sigma^{1/2}\over \pa\nu}
\bigr),&(1.6)\cr}$$
where $\Delta$ is the Laplacian, $\Lambda$ is the voltage-to-current map on
$\pa D$, and $\sigma^{-1/2}$, $\pa\sigma^{1/2}/\pa\nu$ in (1.6) denote the
multiplication operators by the functions $\sigma^{-1/2}\big|_{\pa D}$,
$(\pa\sigma^{1/2}/\pa\nu)\big|_{\pa D}$, respectively. In addition, (1.4)
is always fulfilled if $v$ is given by (1.5) (where $\sigma\ge\sigma_{min}>0$
and $\sigma$ is twice differential in $L^{\infty}(D)$).

Problem 1.1 includes, in particular, the following questions: (a) uniqueness,
(b) reconstruction, (c) stability.

Global uniqueness theorems for Problem 1.1 (in its Calderon or Gel'fand form)
in dimension $d\ge 3$ were obtained for the first time in [23] and [16].
In particular, according to  [16], under assumptions (1.2), (1.4), the map
$\Phi$ uniquely determines $v$.

A global reconstruction method for Problem 1.1 in dimension $d\ge 3$ was
proposed for  the  first time in [16]. In its simplest form, this method of
[16] consists in reducing Problem 1.1 to Problem 4.1 via formulas and
equations  (5.1), (5.2) for $v_2=v$, $\psi_2=\psi$, $h_1=h$, $v_1\equiv 0$,
$\psi_1(x,k)=e^{ikx}$, $R_1(x,y,k)=G(x-y,k)$, $h_1\equiv 0$
  and in solving Problem 4.1 via
formula (4.10), see Sections 3,4 and 5 of the present paper.

Global stability estimates for Problem 1.1 in dimension $d\ge 3$ were given
for the first time in [1]. A variation of this result is presented as
Theorem 2.1 of Section 2 of the present paper.

We recall that in global results one does not assume that $\sigma$ or $v$
is small in some sense or that $\sigma$ or $v$ is
peicewise constant or piecewise real-analytic. For
peicewise constant or piecewise real-analytic  $\sigma$ the first uniqueness
 results for the Calderon version of

\noindent
Problem 1.1 in dimension $d\ge 2$ were given in [6], [11].

As regards global results given in the literature on Problem 1.1 in dimension
$d=2$, the reader is refered to Corollary 2 of [16] and to [15], [3], [13], 
[12], [4] and to references given in [4].

In the present work we continue studies of [16], [17], [19], [21], [20]. One
of the main purposes of these studies was developing an efficient
reconstruction  algorithm for

\noindent
Problem 1.1 in dimension $d=3$. Note that the
aforementioned global reconstruction method of [16] is fine in the sense that
it consists in solving Fredholm linear integral equations of the second type
 and using explicit formulas but this reconstruction is not optimal with
respect to its stability properties. An effectivization of this reconstruction
 of [16] with respect to its stability properties was developed in
[17], [19], [20]. In [21] and in the present work we illustrate our progress
in stability by proving new stability estimates for Problem 1.1 in 3D and,
in particular, by proving Theorem 2.2 of Section 2 in 3D. For sufficiently
regular
potentials this theorem is a principle improvement of the aforementioned
Alessandrini stability result of [1].

Note that algorithms developed in [17], [19], [20] for finding $v$ from
$\Phi$ in 3D can be used even if $v$ has discontinuities. In this case
these algorithms can be considered as methods for finding the smooth part of
$v$ from $\Phi$. In particular, the global algorithms of [17], [20] can be
considered as methods for approximate finding the Fourier
transform $\hat v(p)$, $|p|\le 2\rho$, from $\Phi$ for sufficiently great
$\rho$, where $\hat v(p)$  is defined by (4.12), $d=3$. However, for the case
when $v$ is not smooth enough we did not manage yet to formalize our
principle progress of [20] in global stability as a rigorous mathematical
theorem.

Note that in [21] our new stability estimates of Theorem 2.2 of the next
section were proved in the Born approximation (that is in the linear
approximation near zero potential) only. Besides, a scheme of proof of these
estimates was also mentioned in  [21] for potentials with sufficiently small
norm in dimension $d=3$. In the present work we give a complete proof of
Theorem 2.2 in the general (or by other words global) case in dimension $d=3$.

Our  new stabilty estimates for Problem 1.1 (that is the estimates of Theorem
2.2) are presented and discussed in Section 2. In addition, relations
between Theorems 2.1, 2.2 and global reconstructions of [16], [17], [20] are
explained in Section 3.
The proof of
Theorem 2.2 for the global case in 3D is given in Section 8.

\vskip 2 mm
{\bf 2. Stability estimates}

As in [21] we assume for simplicity that
$$\eqalign{
&D\ \ {\rm is\ an\ open\ bounded\ domain\ in}\ \ \R^d,\ \pa D\in C^2,\cr
&v\in W^{m,1}(\R^d)\ \ {\rm for\ some}\ m>d,\ supp\,v\subset D,\ d\ge 2,\cr}
\eqno(2.1)$$
where
$$W^{m,1}(\R^d)=\{v:\ \pa^Jv\in L^1(\R^d),\ |J|\le m\},\ \ m\in\N\cup 0,
\eqno(2.2)$$
where
$$J\in (\N\cup 0)^d,\ \ |J|=\sum_{i=1}^dJ_i,\ \
\pa^Jv(x)={\pa^{|J|}v(x)\over \pa x_1^{J_1}\ldots \pa x_d^{J_d}}.$$
Let
$$\|v\|_{m,1}=\max_{|J|\le m}\|\pa^Jv\|_{L^1(\R^d)}.\eqno(2.3)$$
Let
$$\eqalign{
&\|A\|\ \ {\rm denote\ the\ norm\ of\ an\ operator}\cr
&A:\ L^{\infty}(\pa D)\to L^{\infty}(\pa D).\cr} \eqno(2.4)$$
We recall that if $v_1$, $v_2$ are potentials satisfying (1.2), (1.3),
where $D$ is fixed, then
$$\Phi_1-\Phi_2\ \ {\rm is\ a\ compact\ operator\ in}\ \ L^{\infty}(\pa D),
\eqno(2.5)$$
where $\Phi_1$, $\Phi_2$ are the DtN maps for $v_1$, $v_2$ respectively,
see [16], [17]. Note also that $(2.1)\Rightarrow (1.2)$.

\vskip 2 mm
{\bf Theorem 2.1} (variation of the result of [1]).
{\it  Let conditions} (1.4), (2.1) {\it hold for potentials} $v_1$
{\it and} $v_2$, {\it where} $D$ {\it is fixed}, $d\ge 3$. {\it Let}
$\|v_j\|_{m,1}\le N$, $j=1,2$, {\it for some} $N>0$. {\it Let} $\Phi_1$,
$\Phi_2$ {\it denote the DtN maps for} $v_1$, $v_2$, {\it respectively. Then}
$$\|v_1-v_2\|_{L^{\infty}(D)}\le
c_1(\ln(3+\|\Phi_1-\Phi_2\|^{-1}))^{-\alpha_1},\eqno(2.6)$$
{\it where} $c_1=c_1(N,D,m)$, $\alpha_1=(m-d)/m$, $\|\Phi_1-\Phi_2\|$
{\it is defined according to} (2.4).

As it was mentioned in [21], Theorem 2.1 follows from formulas (4.9)-(4.11),
(5.1) (of Sections 4 and 5).

A disadvantage of estimate (2.6) is that
$$\alpha_1<1\ \ {\rm for\ any}\ \  m>d\ \ {\rm even\ if}\ \ m\ \ {\rm is\
 very\ great}.\eqno(2.7)$$

\vskip 2 mm
{\bf Theorem 2.2.}
{\it Let the assumptions of Theorem 2.1 hold. Then}
$$\|v_1-v_2\|_{L^{\infty}(D)}\le
c_2(\ln(3+\|\Phi_1-\Phi_2\|^{-1}))^{-\alpha_2},\eqno(2.8)$$
{\it where} $c_2=c_2(N,D,m)$, $\alpha_2=m-d$, $\|\Phi_1-\Phi_2\|$
{\it is defined according to} (2.4).

A principal advantage of estimate (2.8) in comparison with (2.6) is that
$$\alpha_2\to +\infty\ \ {\rm as}\ \  m\to +\infty, \eqno(2.9)$$
in contrast with (2.7).

In the Born approximation, that is in the linear approximation near zero
potential, Theorem 2.2 was proved in [21].

For sufficiently small $N$ in dimension $d=3$,
a scheme of proof of Theorem 2.2 was also  mentioned in [21]. This scheme
involves, in particular, results of [17], [19].

In the general   (or by other words global) case Theorem 2.2 in dimension
$d=3$ is proved in Section 8. This proof involves, in particular, results of
[17], [20].
We see no principal difficulties (except restrictions in time) to give a
similar proof in dimension $d>3$.

We would like to mention that, under the assumptions of Theorems 2.1 and 2.2,
according to the Mandache results of [13], the estimate (2.8) can not hold
with

\noindent
$\alpha_2>m(2d-1)/d$. However, we think that this instability result of
[13] can be sharpened considerably and that in this sense our Theorem 2.2 is
almost optimal already.

For additional information concerning stability and instability results given
in the literature for Problem 1.1 (in its Calderon or Gel'fand form) the
reader is refered to [1], [13], [21], [22], [24] and references therein.
In particular, in [24] estimate (2.6) is proved for $d=2$ and any
$\alpha_1\in ]0,1[$, under the assumptions that
$v_j\in C^2(\bar D)$, $\|v_j\|_{C^2(\bar D)}\le N$ (and, at least, under
the additional assumption that $supp\,v_j\subseteq\bar D_0\subset D$ for some
closed $\bar D_0$), $j=1,2$, where $c_1$ depends on $D$, $N$, $\alpha_1$ (and
on $\bar D_0$).

\vskip 2 mm
{\bf 3. Theorems 2.1 and  2.2 and global reconstructions of [16], [17], [20]}

\vskip 2 mm
In this section we explain relations between Theorems 2.1, 2.2 and global
reconstructions of [16], [17], [20] for Problem 1.1, where we assume that
dimension $d=3$. (The case $d>3$ is similar to the case $d=3$.) The scheme
of the aforementioned global reconstructions consists in the following:
$$\Phi\to\ h\big|_{\Xi}\to\hat v_{\rho,\tau}\to v,\eqno(3.1)$$
where $h$ is the Faddeev generalized scattering amplitude (at zero energy)
defined in the complex domain $\Theta$ (see formulas (4.4), (4.5), (4.1)-(4.3)
 of Section 4 for $d=3$), $\Xi$ is an appropriate subset of $\Theta$, for example
$\Xi$ is the subset of $\Theta$ with the imaginary part length equal to
$\rho$, $\hat v_{\rho,\tau}$ is an approximation to $\hat v$ on
${\cal B}_{2\tau\rho}$,
$$\hat v_{\rho,\tau}\to\hat v\ \ {\rm on}\ \
{\cal B}_{2\tau\rho}\ \ {\rm as}\ \  \rho\to +\infty,\eqno(3.2)$$
where $\hat v$ is the Fourier transform of $v$ (see formula (4.12), $d=3$),
${\cal B}_r$ is the ball of the radius $r$, $\tau$ is a fixed parameter,
$0<\tau\le 1$.

\vskip 2 mm
{\it 3.1. Theorem 2.1 and global reconstructions of [16], [17].}
For the case when
$$\Xi=b\Theta_{\rho}\ \ {\rm (see\ formulas\ (6.2))},\eqno(3.3a)$$
$$\eqalign{
&h\big|_{b\Theta_{\rho}}\ \ {\rm is\ determined\ from}\ \ \Phi\ \ {\rm via}\cr
&(5.1), (5.2)\ \ {\rm considered\ for}\ \ v_2=v,\ \psi_2=\psi,\ h_2=h,\cr
&v_1\equiv 0,\ \psi_1(x,k)=e^{ikx},\ R_1(x,y,k)=G(x-y,k),\ h_1\equiv 0,\cr}
\eqno(3.3b)$$
$$\eqalign{
&\hat v_{\rho,\tau}(p)=h(k(p),l(p)),\ p\in {\cal B}_{2\tau\rho},\ \tau=1,\cr
&{\rm for\ some\ vector-functions}\ \ k(p),\ l(p)\cr
&{\rm such\ that}\ \ k(p),l(p)\in b\Theta_{\rho},\ \ k(p)-l(p)=p,\cr}
\eqno(3.3c)$$
reconstruction (3.1) was, actually, given in [16]. In this case (3.2) holds
according to (4.10), (4.11). (In the framework of [16] less precise version
of (4.11) was known.)

General formulas and equations (5.1)-(5.3) for finding
$h\big|_{b\Theta_{\rho}}$ from $\Phi$ with known background potential
$v_1\not\equiv 0$ were given for the first time in [17]. If unknown potential
$v$ is sufficiently close to some known background potential
$v_1\not\equiv 0$, then formulas and equations (5.1)-(5.3) considered for
$v_2=v$, $\psi_2=\psi$, $h_2=h$ give more stable determination of
$h\big|_{b\Theta_{\rho}}$ from $\Phi$  than these formulas and equations for
$v_1\equiv 0$ (see [17] for more detailed discussion of this issue).

Note that Theorem 2.1 follows from formulas (4.10), (4.11), (5.1) mentioned
above in connection with the versions of (3.1), (3.2) of [16], [17].

\vskip 2 mm
{\it 3.2. Theorem 2.2 and global reconstructions of [17], [20].}
The main disadvantage of the global reconstruction (3.1), (3.2) in the
framework of [16], [17] is related with the following two facts:

(1) The determination of $h\big|_{\Xi}$ for $\Xi=b\Theta_{\rho}$ or, more
generally, for $\Xi\subseteq\bar\Theta_{\rho}$,
$\Xi\cap b\Theta_{\rho}\ne\emptyset$ (see formulas (6.2)) from $\Phi$ via
formulas and equations (5.1)-(5.3) for $v_2=v$, $\psi_2=\psi$,
$h_2=h$ with some known background potential $v_1$
is stable for
relatively small $\rho$, but is very unstable for $\rho\to +\infty$ in
the points of $\Xi$ with sufficiently great imaginary part (see [17], [19],
[20] for more detailed discussion of this issue).

(2) For $\hat v_{\rho,\tau}$ defined as in (3.3c) the decay of the error
$\hat v-\hat v_{\rho,\tau}$ on ${\cal B}_{2\rho\tau}$ is very slow for
$\rho\to +\infty$ (not faster than $O(\rho^{-1})$ even for infinitely
smooth compactly supported $v$) (see [20] for more detailed discussion of
this issue).

As a corollary, the global reconstruction (3.1), (3.2) in the framework of
[16] and even in the framework of [17] is not optimal with respect to its
stability properties.
A principle effectivization of the global reconstruction (3.1), (3.2) with
respect to its stability properties was given in [20].

The main new result of [20] consists in stable construction of some
$\hat v_{\rho,\tau}$ on ${\cal B}_{2\tau\rho}$ from $h$ on
$b\Theta_{\rho,\tau}\subseteq b\Theta_{\rho}$ (see definitions (6.2)) in
such a way that
$$\eqalign{
&|\hat v(p)-\hat v_{\rho,\tau}(p)|\le C(m,\mu_0,\tau)N^2(1+|p|)^{-\mu_0}
\rho^{-(m-\mu_0)}\cr
&{\rm for}\ \ p\in {\cal B}_{2\tau\rho}, \ \rho\ge\rho_1\ \ {\rm and\ fixed}\
\tau\in ]0,\tau_1[,\cr}\eqno(3.4)$$
under the assumptions that
$$\eqalign{
&v\in W^{m,1}(\R^2),\ \ \|v\|_{m,1}<N\ \ {\rm (see\ definitions}\ \
(2.2), (2.3)),\cr
&m>2,\ 2\le\mu_0<m,\ \tau_1=\tau_1(m,\mu_0,N),\ \ \rho_1=\rho_1(m,\mu_0,N),
\cr}\eqno(3.5)$$
where $C(m,\mu_0,\tau)$,
$\tau_1(m,\mu_0,N)$, $\rho_1(m,\mu_0,N)$ are special constants and,
in particular, $0<\tau_1<1$, see Theorem 2.1 of [20]. The scheme of
this construction of [20] consists in the following:
$$\eqalign{
&h\big|_{b\Theta_{\rho,\tau}}\to H\big|_{b\Omega_{\rho,\tau}}\to
H\big|_{b\Lambda_{\rho,\tau,\nu}}\cr
&\to  H^0\big|_{\Lambda_{\rho,\tau,\nu}}\to\tilde H_{\rho,\tau}\ \ {\rm on}\
\Lambda_{\rho,\tau,\nu}\cr
&\to \hat v_{\rho,\tau}\ \ {\rm on}\ \ {\cal B}_{2\tau\rho}.\cr}\eqno(3.6)$$
Here $h\big|_{b\Theta_{\rho,\tau}}$ and $H\big|_{b\Omega_{\rho,\tau}}$ are
related by (4.13) (see also the definitions of $b\Theta_{\rho,\tau}$ and
$b\Omega_{\rho,\tau}$ given in (6.2), (6.3)),
$H\big|_{b\Lambda_{\rho,\tau,\nu}}$ is defined in terms of
$H\big|_{b\Omega_{\rho,\tau}}$ by means of (7.4) (where
$b\Lambda_{\rho,\tau,\nu}$ is defined in (6.13)),
$H^0\big|_{\Lambda_{\rho,\tau,\nu}}$ is defined in terms of
$H\big|_{b\Lambda_{\rho,\tau,\nu}}$ by (7.8) (see also formulas (6.13),
(6.16)), $\tilde H_{\rho,\tau}$ is defined as the solution of (7.7) for
$H$ on  $\Lambda_{\rho,\tau,\nu}$, where $Q_{\rho,\tau}$ is replaced by zero
in (7.7), finally, $\hat v_{\rho,\tau}$ on  ${\cal B}_{2\tau\rho}$ is
defined as $\hat v_{\rho,\tau}=\hat v^+_{\rho,\tau}$ or as
$\hat v_{\rho,\tau}=\hat v^-_{\rho,\tau}$, where
$$\eqalign{
&\hat v^+_{\rho,\tau}(p)=\lim\limits_{\lambda\to 0}
\tilde H_{\rho,\tau}(\lambda,p),\cr
&\hat v^-_{\rho,\tau}(p)=\lim\limits_{\lambda\to \infty}
\tilde H_{\rho,\tau}(\lambda,p),\cr}\eqno(3.7)$$
where definitions (3.7) are similar to formulas (7.16). For more detailed
discussion of (3.6) see Section 5 of [20].

Note that in [20] we deal with approximate finding $v$ from
$h\big|_{b\Theta_{\rho}}$ in a similar way with approximate finding
$v$ from the scattering amplitude $f$ at fixed energy $E$ of [18].
The parameter $\rho$ of [20] plays the role of the parameter $\sqrt{E}$
of [18].

One can see that for sufficiently regular potentials $v$ that is for
$v\in W^{m,1}(\R^3)$, where $m$ is sufficiently great, and for
$\hat v_{\rho,\tau}$ constructed from $h\big|_{b\Theta_{\rho,\tau}}$ via (3.6)
 as in
[20] the decay rate of the error $\hat v-\hat v_{\rho,\tau}$ on
${\cal B}_{2\tau\rho}$ is $O(\rho^{-(m-2)})$ as $\rho\to +\infty$ (see (3.4)
for $\mu_0=2$) and is much faster than for  $\hat v_{\rho,\tau}$
given by (3.3c) as in [16], [17]. Therefore, for Problem 1.1 in dimension $d=3$ the global
reconstruction of [20] (using also (5.1)-(5.3) for finding
$h\big|_{b\Theta_{\rho,\tau}}$ from $\Phi$ as in [17]) has much more optimal
stability
properties than the global reconstructions of [16], [17].

Therefore, Theorem 1.2 illustrating the progress in global stability of
[17], [20] is much more optimal than Theorem 1.1 related with global
reconstruction results of [16] and [17] only.

\vskip 2 mm
{\bf 4. Faddeev functions}

We consider the Faddeev functions $G$, $\psi$ and $h$ (see [7], [8],
[10], [16]):
$$\eqalignno{
&\psi(x,k)=e^{ikx}+\int_{\R^d}G(x-y,k)v(y)\psi(y,k)dy,&(4.1)\cr
&G(x,k)=e^{ikx}g(x,k),\ \
g(x,k)=-(2\pi)^{-d}\int_{\R^d}{e^{i\xi x}d\xi\over
{\xi^2+2k\xi}},&(4.2)\cr}$$
where $x\in\R^d$, $k\in\Sigma$,
$$\eqalignno{
&\Sigma=\{k\in\C^d:\ \ k^2=k_1^2+\ldots +k_d^2=0\};&(4.3)\cr
&h(k,l)=(2\pi)^{-d}\int_{\R^d}e^{-ilx}v(x)\psi(x,k)dx,&(4.4)
\cr}$$
where $(k,l)\in\Theta$,
$$\Theta=\{k\in\Sigma,\ \ l\in\Sigma:\ Im\,k=Im\,l\}.\eqno(4.5)$$
One can consider (4.1), (4.4) assuming that $v$ is a sufficiently regular
function on $\R^d$ with sufficient decay at infinity. For example, one can
consider (4.1), (4.4) assuming that (1.2) holds.

We recall that:
$$\Delta G(x,k)=\delta(x),\ \ x\in\R^d,\ \ k\in\Sigma;\eqno(4.6)$$
formula (4.1) at fixed $k$ is considered as an equation for
$$\psi=e^{ikx}\mu(x,k),\eqno(4.7)$$
where $\mu$ is sought in $L^{\infty}(\R^d)$; as a corollary of (4.1),(4.2),
(4.6), $\psi$
satisfies (1.1); $h$ of (4.4) is a generalized "scattering" amplitude
in the complex domain at zero energy.

Note that, actually, $G$, $\psi$, $h$ of (4.1)-(4.5) are zero energy
restrictions of functions introduced by Faddeev as extentions to the
complex domain of some functions of the classical scattering theory
for the Schr\"odinger equation at positive energies. In addition,
$G$, $\psi$, $h$ in their zero energy restriction were considered for the
first time in [2]. The Faddeev functions $G$, $\psi$, $h$ were, actually,
rediscovered in [2].

We recall also that, under the assumptions of Theorem 2.1,
$$\mu(x,k)\to 1\ \ {\rm as}\ \ |Im\,k|\to\infty\ \ ({\rm uniformly\ in}\ \ x)
\eqno(4.8)$$
and, for any $\sigma>1$,
$$|\mu(x,k)|<\sigma \ \ {\rm for}\ \ |Im\,k|\ge r_1(N,D,m,\sigma),\eqno(4.9)$$
where $x\in\R^d$, $k\in\Sigma$;
$$\hat v(p)=\lim_{\scriptstyle (k,l)\in\Theta,\ k-l=p\atop\scriptstyle
|Im\,k|=|Im\,l|\to\infty}h(k,l)\ \ {\rm for\ any}\ \ p\in\R^d,\eqno(4.10)$$
$$\eqalign{
&|\hat v(p)-h(k,l)|\le {c_3(D,m)N^2\over \rho}\ \ {\rm for}\ \
(k,l)\in\Theta,\ p=k-l,\cr
&|Im\,k|=|Im\,l|=\rho\ge r_2(N,D,m),\cr}\eqno(4.11)$$
where
$$\hat v(p)=\bigl({1\over 2\pi}\bigr)^d\int_{\R^d}e^{ipx}v(x)dx,\ \ p\in\R^d.
\eqno(4.12)$$

Results of the type (4.8), (4.9) go back to [2]. Results of the type
(4.10), (4.11) (with less precise right-hand side in (4.11)) go back to [10].
Estimates (4.8), (4.11) are  related also with some important
$L_2$-estimate going back to [23] on the Green function $g$ of (4.1).

Note also that in some considerations it is convenient to consider $h$     on
$\Theta$ as $H$ on $\Omega$, where
$$\eqalign{
&h(k,l)=H(k,k-l),\ \ (k,l)\in\Theta,\cr
&H(k,p)=h(k,k-p),\ \ (k,p)\in\Omega,\cr}\eqno(4.13)$$
$$\Omega=\{k\in\C^d,\ p\in\R^d:\ k^2=0,\ p^2=2kp\}.\eqno(4.14)$$

For more information on properties of the Faddeev functions
$G$, $\psi$, $h$, see [10], [17], [20] and references therein.

In the next section we recall that Problem 1.1 (of Introduction) admits a
reduction to the following inverse "scattering" problem:

\vskip 2 mm
{\bf Problem 4.1.}
Given $h$ on $\Theta$, find $v$ on $\R^d$.

\vskip 2 mm
{\bf 5. Reduction of [16], [17] of Problem 1.1 to Problem 4.1}

Let conditions (1.2), (1.4) hold for potentials $v_1$ and $v_2$, where
$D$ is fixed. Let $\Phi_i$, $\psi_i$, $h_i$ denote the DtN map $\Phi$
and the Faddeev functions $\psi$, $h$ for $v=v_i$, $i=1,2$. Let also
$\Phi_i(x,y)$ denote the Schwartz kernel $\Phi(x,y)$ of the integral operator
$\Phi$ for $v=v_i$, $i=1,2$. Then (see [17] for details):
$$h_2(k,l)-h_1(k,l)=\bigl({1\over 2\pi}\bigr)^d
\int\limits_{\pa D}\!\!\int\limits_{\pa D}
\psi_1(x,-l)(\Phi_2-\Phi_1)(x,y)\psi_2(y,k)dydx,\eqno(5.1)$$
where $(k,l)\in\Theta$;
$$\eqalignno{
&\psi_2(x,k)=\psi_1(x,k)+\int_{\pa D}A(x,y,k)\psi_2(y,k)dy,\ x\in\pa D,
&(5.2a)\cr
&A(x,y,k)=\int_{\pa D}R_1(x,z,k)(\Phi_2-\Phi_1)(z,y)dz,\ x,y\in\pa D,
&(5.2b)\cr
&R_1(x,y,k)=G(x-y,k)+\int_{\R^d}G(x-z,k)v_1(z)
R_1(z,y,k)dz,\ x,y\in\R^d,
&(5.3)\cr}$$
where $k\in\Sigma$. Note that: (5.1) is an explicit formula, (5.2a) is
considered as an equation for finding $\psi_2$ on $\pa D$ from $\psi_1$ on
$\pa D$ and $A$ on $\pa D\times\pa D$ for each fixed $k$, (5.2b)
is an explicit formula, (5.3) is an equation for finding $R_1$ from $G$
and $v_1$, where $G$ is the function of (4.2).

Note that formulas and equations (5.1)-(5.3) for $v_1\equiv 0$ were given
in [16] (see also [10] (Note added in proof), [14], [15]).
In this case $h_1\equiv 0$, $\psi_1=e^{ikx}$, $R_1=G(x-y,k)$. Formulas and
equations (5.1)-(5.3) for the general case were given in [17].

Formulas and equations (5.1)-(5.3) with fixed background potential $v_1$
reduce Problem 1.1 (of Introduction) to Problem 4.1 (of Section 4).
In this connection we consider (5.1)-(5.3) for $v_2=v$, $\psi_2=\psi$,
$h_2=h$ and this reduction consists of the following steps:
\item{(1)} $v_1\to \Phi_1$, $\psi_1$, $R_1$, $h_1$ via formulas and
equations of the diret problem for $v_1$;
\item{(2)} $\Phi$, $\Phi_1$, $R_1\big|_{\pa D\times\pa D\times\Sigma}\to A$
via (5.2b);
\item{(3)} $A$,
$\psi_1\big|_{\pa D\times\Sigma}\to \psi\big|_{\pa D\times\Sigma}$ via (5.2a);
\item{(4)} $h_1$, $\psi_1\big|_{\pa D\times\Sigma}$,
$\psi\big|_{\pa D\times\Sigma}$, $\Phi$, $\Phi_1\to h$ via (5.1).

Here, on the step (4) we find $h$ on $\Theta$ for the unknown potential
$v$ of Problem 1.1.

\vskip 2 mm
{\bf 6. Some considerations related with $\Theta$ and $\Omega$}

\vskip 2 mm
{\it 6.1 Some subsets of $\Theta$ and $\Omega$.}
Let
$$\eqalign{
&{\cal B}_r=\{p\in\R^d:\ |p|<r\},\ \ \pa {\cal B}_r=\{p\in\R^d:\ |p|=r\},\cr
&\bar {\cal B}_r={\cal B}_r\cup\pa {\cal B}_r,\ \ {\rm where}\ \ r>0.\cr}
\eqno(6.1)$$
In addition to $\Theta$ of (4.5), we consider, in particular, the following
its subsets:
$$\eqalign{
&\Theta_{\rho}=\{(k,l)\in\Theta:\ |Im\,k|=|Im\,l|<\rho\},\cr
&b\Theta_{\rho}=\{(k,l)\in\Theta:\ |Im\,k|=|Im\,l|=\rho\},\cr
&\bar\Theta_{\rho}=\Theta_{\rho}\cup b\Theta_{\rho},\cr
&\Theta_{\rho,\tau}^{\infty}=\{(k,l)\in\Theta\b\bar\Theta_{\rho}:\
k-l\in {\cal B}_{2\rho\tau}\},\cr
&b\Theta_{\rho,\tau}=\{(k,l)\in b\Theta_{\rho}:\
k-l\in {\cal B}_{2\rho\tau}\},\cr}\eqno(6.2)$$
where $\rho>0$, $0<\tau<1$, and ${\cal B}_r$ is defined in (6.1).

In addition to $\Omega$ of (4.14), we consider, in particular, the following
its subsets:
$$\eqalign{
&\Omega_{\rho}=\{(k,p)\in\Omega:\ |Im\,k|<\rho\},\cr
&b\Omega_{\rho}=\{(k,p)\in\Omega:\ |Im\,k|=\rho\},\cr
&\bar\Omega_{\rho}=\Omega_{\rho}\cup b\Omega_{\rho},\cr
&\Omega_{\rho,\tau}^{\infty}=\{(k,p)\in\Omega\b\bar\Omega_{\rho}:\
p\in {\cal B}_{2\rho\tau}\},\cr
&b\Omega_{\rho,\tau}=\{(k,p)\in b\Omega_{\rho}:\
p\in {\cal B}_{2\rho\tau}\},\cr}\eqno(6.3)$$
where $\rho>0$, $0<\tau<1$, and ${\cal B}_r$ is defined in (6.1).

Note that
$$\eqalign{
&\Omega\approx\Theta,\ \Omega_{\rho}\approx\Theta_{\rho},\
b\Omega_{\rho}\approx b\Theta_{\rho},\cr
&\Omega_{\rho,\tau}^{\infty}\approx\Theta_{\rho,\tau}^{\infty},\
b\Omega_{\rho,\tau}\approx b\Theta_{\rho,\tau},\cr}\eqno(6.4)$$
or, more precisely,
$$\eqalign{
&(k,p)\in\Omega\Rightarrow (k,k-p)\in\Theta,\ \
(k,l)\in\Theta\Rightarrow (k,k-l)\in\Omega\cr
&{\rm and\ the\ same\ for}\ \ \Omega_{\rho},\ b\Omega_{\rho},\
\Omega_{\rho,\tau}^{\infty},\ b\Omega_{\rho,\tau}\cr
&{\rm and}\ \ \Theta_{\rho},\ b\Theta_{\rho},\
\Theta_{\rho,\tau}^{\infty},\ b\Theta_{\rho,\tau},\ \ {\rm respectively,\ in\
place\ of}\ \ \Omega\ \ {\rm and}\ \ \Theta.\cr}\eqno(6.5)$$
We consider also, in particular,
$$\eqalign{
&\Omega_{\nu}=\{(k,p)\in\Omega:\ p\notin {\cal L}_{\nu}\},\cr
&\Omega_{\rho,\tau,\nu}^{\infty}=\Omega_{\rho,\tau}^{\infty}\cap\Omega_{\nu},
\ b\Omega_{\rho,\tau,\nu}=b\Omega_{\rho,\tau}\cap\Omega_{\nu},\cr}\eqno(6.6)$$
where
$${\cal L}_{\nu}=\{p\in\R^d:\ p=t\nu,\ t\in\R\},\eqno(6.7)$$
$\nu\in\S^{d-1}$, $\rho>0$, $0<\tau<1$.

\vskip 2 mm
{\it 6.2. Coordinates on $\Omega$ for $d=3$.}
In this subsection we assume that $d=3$ in formulas (4.5), (4.14),
(6.1)-(6.7).

For $p\in\R^3\b {\cal L}_{\nu}$ we consider $\theta(p)$ and $\omega(p)$ such
that
$$\eqalign{
&\theta(p),\omega(p)\ \ {\rm smoothly\ depend\ on}\ \ p\in\R^3\b
{\cal L}_{\nu},\cr
&{\rm take\ values\ in}\ \ \S^2,\ \ {\rm and}\cr
&\theta(p)p=0,\ \omega(p)p=0,\ \theta(p)\omega(p)=0,\cr}\eqno(6.8)$$
where ${\cal L}_{\nu}$ is defined by (6.7) (for $d=3$).

Assumptions (6.8) imply that
$$\omega(p)={p\times\theta(p)\over |p|}\ \ {\rm for}\ \
p\in\R^3\b {\cal L}_{\nu} \eqno(6.9a)$$
or
$$\omega(p)=-{p\times\theta(p)\over |p|}\ \ {\rm for}\ \
p\in\R^3\b {\cal L}_{\nu}, \eqno(6.9b)$$
where $\times$ denotes vector product.

To satisfy (6.8), (6.9a) we can take
$$\theta(p)={\nu\times p\over |\nu\times p|},\ \omega(p)=
{p\times\theta(p)\over |p|},\ p\in\R^3\b {\cal L}_{\nu}.\eqno(6.10)$$

Let $\theta,\omega$ satisfy (6.8). Then (according to [19]) the following
formulas give a diffeomorphism between $\Omega_{\nu}$ and
$(\C\b 0)\times (\R^3\b {\cal L}_{\nu})$:
$$(k,p)\to (\lambda,p),\ \ {\rm where}\ \
\lambda=\lambda(k,p)={2k(\theta(p)+i\omega(p))\over i|p|},\eqno(6.11a)$$
$$\eqalign{
&(\lambda,p)\to (k,p),\ \ {\rm where}\ \ k=k(\lambda,p)=\kappa_1(\lambda,p)
\theta(p)+\kappa_2(\lambda,p)\omega(p)+{p\over 2},\cr
&\kappa_1(\lambda,p)={i|p|\over 4}(\lambda+{1\over \lambda}),
\ \  \kappa_2(\lambda,p)={|p|\over 4}(\lambda-{1\over \lambda}),\cr}
\eqno(6.11b)$$
where $(k,p)\in\Omega_{\nu}$,
$(\lambda,p)\in (\C\b 0)\times (\R^3\b {\cal L}_{\nu})$.
In addition, formulas (6.11a), (6.11b) for $\lambda(k)$ and $k(\lambda)$
at fixed $p\in\R^3\b {\cal L}_{\nu}$
give a diffeomorphism between $Z_p=\{k\in\C^3:\ (k,p)\in\Omega\}$ for fixed
$p$ and $\C\b 0$.

In addition, for $k$ and $\lambda$ of (6.11) we have that
$$|Im\,k|={|p|\over 4}\bigl(|\lambda|+{1\over |\lambda|}\bigr),\ \
|Re\,k|={|p|\over 4}\bigl(|\lambda|+{1\over |\lambda|}\bigr),\eqno(6.12)$$
where $(k,p)\in\Omega_{\nu}$,
$(\lambda,p)\in (\C\b 0)\times (\R^3\b {\cal L}_{\nu})$.

Let
$$\Lambda_{\rho,\nu}=\{(\lambda,p):\ \lambda\in {\cal D}_{\rho/|p|},\ \
p\in\R^3\b {\cal L}_{\nu}\},\eqno(6.13)$$
$$\eqalign{
&\Lambda_{\rho,\tau,\nu}=\{(\lambda,p):\ \lambda\in {\cal D}_{\rho/|p|},\ \
p\in\R^3\b {\cal L}_{\nu},\ \ |p|<2\tau\rho\},\cr
&b\Lambda_{\rho,\tau,\nu}=\{(\lambda,p):\ \lambda\in {\cal T}_{\rho/|p|},\ \
p\in\R^3\b {\cal L}_{\nu},\ \ |p|<2\tau\rho\},\cr}$$
where $\rho>0$, $0<\tau<1$, $\nu\in\S^2$,
$${\cal D}_r=\{\lambda\in\C\b 0:\ {1\over 4}(|\lambda|+|\lambda|^{-1})>r\},\
r>0,\eqno(6.14)$$
$${\cal T}_r=\{\lambda\in\C:\ {1\over 4}(|\lambda|+|\lambda|^{-1})=r\},\
r\ge 1/2,\eqno(6.15)$$
${\cal L}_{\nu}$ is defined by (6.7) (for $d=3$).

Note that
$$\eqalign{
&\Lambda_{\rho,\tau,\nu}=
\Lambda^+_{\rho,\tau,\nu}\cup \Lambda^-_{\rho,\tau,\nu},\
\Lambda^+_{\rho,\tau,\nu}\cap \Lambda^-_{\rho,\tau,\nu}=\emptyset,\cr
&b\Lambda_{\rho,\tau,\nu}=
b\Lambda^+_{\rho,\tau,\nu}\cup b\Lambda^-_{\rho,\tau,\nu},\cr}\eqno(6.16)$$
where
$$\eqalign{
&\Lambda^{\pm}_{\rho,\tau,\nu}=\{(\lambda,p):\ \
\lambda\in {\cal D}^{\pm}_{\rho/|p|},\ \
p\in {\cal B}_{2\tau\rho}\b {\cal L}_{\nu}\},\cr
&b\Lambda^{\pm}_{\rho,\tau,\nu}=\{(\lambda,p):\ \
\lambda\in {\cal T}^{\pm}_{\rho/|p|},\ \
p\in {\cal B}_{2\tau\rho}\b {\cal L}_{\nu}\},\cr}\eqno(6.17)$$
$$\eqalign{
&{\cal D}^{\pm}_r=\{\lambda\in\C\b 0:\ \
{1\over 4}(|\lambda|+|\lambda|^{-1})>r,\ \ |\lambda|^{\pm 1}<1\},\cr
&{\cal T}^{\pm}_r=\{\lambda\in\C:\ \
{1\over 4}(|\lambda|+|\lambda|^{-1})=r,\ \ |\lambda|^{\pm 1}\le 1\},\ r>1/2,
\cr}\eqno(6.18)$$
where $\rho>0$, $\tau\in ]0,1[$, $\nu\in\S^2$.

Using (6.12) one can see that  formulas (6.11) give also the following
diffeomorphisms
$$\eqalign{
&\Omega_{\nu}\b\bar\Omega_{\rho}\approx\Lambda_{\rho,\nu},\ \
\Omega_{\rho,\tau,\nu}^{\infty}\approx\Lambda_{\rho,\tau,\nu},\cr
&b\Omega_{\rho,\tau,\nu}\approx b\Lambda_{\rho,\tau,\nu},\cr
&Z_{p,\rho}^{\infty}=\{k\in\C^3:\ (k,p)\in\Omega_{\nu}\b\bar\Omega_{\rho}\}
\approx {\cal D}_{\rho/|p|}\ \ {\rm for\ fixed}\ \ p,\cr}\eqno(6.19)$$
where $\rho>0$, $0<\tau<1$, $\nu\in\S^2$ (and where we use the definitions
(6.3), (6.6), (6.13)).

In [19] $\lambda,p$ of (6.11) were used as coordinates on
$\Omega$. In the present work we use them also as coordinates on
$\Omega\b\Omega_{\rho}$ (or more precisely on $\Omega_{\nu}\b\Omega_{\rho}$).

\vskip 2 mm
{\bf 7. An integral equation of [20] and some related formulas}

In the main considerations of [20] it is assumed that $d=3$ and the basic
assumption on $v$ consists in the following condition on its Fourier
transform:
$$\hat v\in L_{\mu}^{\infty}(\R^3)\cap {\cal C}(\R^3)\ \ {\rm for\ some\
real}\ \ \mu\ge 2,\eqno(7.1)$$
where $\hat v$ is defined by (4.12) (for $d=3$),
$$\eqalign{
&L_{\mu}^{\infty}(\R^d)=\{u\in L^{\infty}(\R^d):\ \ \|u\|_{\mu}<+\infty\},\cr
&\|u\|_{\mu}=ess \sup\limits_{p\in\R^d}(1+|p|)^{\mu}|u(p)|,\ \ \mu>0,\cr}
\eqno(7.2)$$
and ${\cal C}$ denotes the space of continuous functions.

Note that
$$\eqalign{
&v\in W^{m,1}(\R^d)\Longrightarrow
\hat v\in L_{\mu}^{\infty}(\R^d)\cap {\cal C}(\R^d),\cr
&\|\hat v\|_{\mu}\le c_4(m,d)\|v\|_{m,1}\ \ {\rm for}\ \ \mu=m,\cr}\eqno(7.3)
$$
where $W^{m,1}$, $L_{\mu}^{\infty}$ are the spaces of (2.2), (7.2).

Let
$$H(\lambda,p)=H(k(\lambda,p),p),\ \ (\lambda,p)\in (\C\b 0)\times
(\R^3\b {\cal L}_{\mu}),\eqno(7.4)$$
where $H$ is the function of (4.13), $\lambda$, $p$ are the coordinates of
Subsection 6.2 under assumption (6.9a).

Let
$$\eqalign{
&L_{\mu}^{\infty}(\Lambda_{\rho,\tau,\nu})=
\{U\in L^{\infty}(\Lambda_{\rho,\tau,\nu}):\ \ |||U|||_{\rho,\tau,\mu}<\infty
\},\cr
&|||U|||_{\rho,\tau,\mu}=
ess \sup\limits_{(\lambda,p)\in\Lambda_{\rho,\tau,\nu}}(1+|p|)^{\mu}
|U(\lambda,p)|,\ \ \mu>0,\cr}\eqno(7.5)$$
where $\Lambda_{\rho,\tau,\nu}$ is defined in (6.13), $\rho>0$,
$\tau\in ]0,1[$, $\nu\in\S^2$, $\mu>0$.

Let $v$ satisfy (7.1) and $\|\hat v\|_{\mu}\le C$. Let
$$\eta(C,\rho,\mu)\buildrel \rm def \over = a(\mu)C (\ln\rho)^2\rho^{-1}<1,\
\ln\rho\ge 2,\eqno(7.6)$$
where $a(\mu)$ is the constant $c_2(\mu)$ of [20]. Let $H(\lambda,p)$ be
defined by (7.4) and be considered as a function on
$\Lambda_{\rho,\tau,\nu}$ of (6.13). Then (see Section 5 of [20]):
$$H=H^0+M_{\rho,\tau}(H)+Q_{\rho,\tau},\ \ \tau\in ]0,1[,\eqno(7.7)$$
where
$$\eqalignno{
&H^0(\lambda,p)={1\over 2\pi i}\int\limits_{{\cal T}^+_{\rho/|p|}}
H(\zeta,p){d\zeta\over {\zeta-\lambda}},\
(\lambda,p)\in\Lambda^+_{\rho,\tau,\nu},&(7.8a)\cr
&H^0(\lambda,p)=-{1\over 2\pi i}\int\limits_{{\cal T}^-_{\rho/|p|}}
H(\zeta,p){\lambda d\zeta\over \zeta(\zeta-\lambda)},\
(\lambda,p)\in\Lambda^-_{\rho,\tau,\nu},&(7.8b)\cr}$$
where  $\Lambda^{\pm}_{\rho,\tau,\nu}$, ${\cal T}_r^{\pm}$ are defined
in (6.17), (6.18) (and where the integrals along ${\cal T}_r^{\pm}$ are
taken in the counter-clock wise direction);
$$\eqalign{
&M_{\rho,\tau}(U)(\lambda,p)=M^+_{\rho,\tau}(U)(\lambda,p)=\cr
&-{1\over \pi}\int\!\!\!\int\limits_{{\cal D}^+_{\rho/|p|}}
(U,U)_{\rho,\tau}(\zeta,p){d Re\,\zeta d Im\,\zeta\over {\zeta-\lambda}},\
(\lambda,p)\in\Lambda^+_{\rho,\tau,\nu},\cr}\eqno(7.9a)$$
$$\eqalign{
&M_{\rho,\tau}(U)(\lambda,p)=M^-_{\rho,\tau}(U)(\lambda,p)=\cr
&-{1\over \pi}\int\!\!\!\int\limits_{{\cal D}^-_{\rho/|p|}}
(U,U)_{\rho,\tau}(\zeta,p){\lambda d Re\,\zeta d Im\,\zeta
\over \zeta(\zeta-\lambda)},\
(\lambda,p)\in\Lambda^-_{\rho,\tau,\nu},\cr}\eqno(7.9b)$$
$$\eqalign{
&(U_1,U_2)_{\rho,\tau}(\zeta,p)=\{\chi_{2\tau\rho}U_1^{\prime},
\chi_{2\tau\rho}U_2^{\prime}\}(\zeta,p),\  (\zeta,p)
\in\Lambda_{\rho,\tau,\nu},\cr
&\chi_{2\tau\rho}U_j^{\prime}(k,p)=U_j(\lambda(k,p),p),\
(k,p)\in\Omega^{\infty}_{\rho,\tau,\nu},\cr
&\chi_{2\tau\rho}U_j^{\prime}(k,p)=0,\ \ |p|\ge 2\tau\rho,\ \ j=1,2,\cr}
\eqno(7.10)$$
where $U, U_1, U_2$ are test functions on $\Lambda_{\rho,\tau,\nu}$,
$\Omega^{\infty}_{\rho,\tau,\nu}$ is defined in (6.6), $\lambda(k,p)$ is
defined in (6.11a), $\{\cdot,\cdot\}$ is defined by the formula
$$\eqalign{
&\{F_1,F_2\}(\lambda,p)=-{\pi\over 4}\int\limits_{-\pi}^{\pi}
\bigr({|p|\over 2}{{|\lambda|^2-1}\over \bar\lambda|\lambda|}(\cos\v-1)-
{|p|\over \bar\lambda}\sin\v\bigr)\times\cr
&F_1(k(\lambda,p),-\xi(\lambda,p,\v))
F_2(k(\lambda,p)+\xi(\lambda,p,\v),p+\xi(\lambda,p,\v))d\v,\cr}\eqno(7.11)
$$
for $(\lambda,p)\in\Lambda_{\rho,\nu}$, where $F_1, F_2$ are test functions
on $\Omega\b\bar\Omega_{\rho}$, $k(\lambda,p)$ is defined in (6.11b),
$\Lambda_{\rho,\nu}$ is defined in (6.13),
$$\eqalignno{
&\xi(\lambda,p,\v)=Re\,k(\lambda,p)(\cos\v-1)+k^{\perp}(\lambda,p)\sin\v,
&(7.12)\cr
&k^{\perp}(\lambda,p)={Im\,k(\lambda,p)\times Re\,k(\lambda,p)\over
|Im\,k(\lambda,p)|},&(7.13)\cr}$$
where $\times$ in (7.13) denotes vector product;
$$\eqalignno{
&H, H^0, Q_{\rho,\tau}\in L_{\mu}^{\infty}(\Lambda_{\rho,\tau,\nu}),&(7.14)\cr
&|||H|||_{\rho,\tau,\mu_0}\le {C\over {1-\eta(C,\rho,\mu)}},&(7.15a)\cr
&|||H^0|||_{\rho,\tau,\mu_0}\le {C\over {1-\eta(C,\rho,\mu)}}
\bigl(1+{c_5(\mu_0)C\over {1-\eta(C,\rho,\mu)}}\bigr),&(7.15b)\cr
&|||Q_{\rho,\tau}|||_{\rho,\tau,\mu_0}\le
{3c_5(\mu_0)C^2\over (1-\eta(C,\rho,\mu))^2(1+2\tau\rho)^{\mu-\mu_0}},&(7.15c)
\cr}$$
where $2\le\mu_0\le\mu$, $c_5$ is the constant $b_4$ of [20],
$\eta(C,\rho,\mu)$ is defined by (7.6).

Following [20] we consider (7.7) as an approxiate integral equation for
finding $H$ on  $\Lambda_{\rho,\tau,\nu}$ from $H^0$ on
$\Lambda_{\rho,\tau,\nu}$ with unknown remainder $Q_{\rho,\tau}$.

Note also that if $\hat v$ satisfies (7.1), then (see [19], [20])
$$\eqalign{
&H(\lambda,p)\to\hat v(p)\ \ {\rm as}\ \ \lambda\to 0,\cr
&H(\lambda,p)\to\hat v(p)\ \ {\rm as}\ \ \lambda\to\infty,\cr}\eqno(7.16)$$
where $p\in {\cal B}_{2\tau\rho}\b {\cal L}_{\nu}$, $H$ is defined by (7.4)
and is considered as a function on $\Lambda_{\rho,\tau,\nu}$, $\rho>0$,
$0<\tau<1$, $\nu\in\S^2$.

\vskip 2 mm
{\bf 8. Proof of Theorem 2.2 for $d=3$}

Our proof of Theorem 2.2 for $d=3$ is based on formulas (7.16) of Section 7
and on Lemma 8.4 given below in the present section. In turn, Lemma 8.4
follows from formula (7.3) of Section 7 and from Lemmas 8.1, 8.2, 8.3 given
below in the present section. In turn:
\item{(1)} Lemma 8.1 follows from some estimates and Lemmas of [20]
related with studies of the non-linear integral equation given as equation
 (7.7) of the present paper;
\item{(2)} Lemma 8.2 does not follow immediately from results of preceding
works and its proof is given in Section 9;
\item{(3)} Lemma 8.3 follows from formulas (4.9), (5.1) and from definitions,
see Section 9 for proof of this lemma.

\vskip 2 mm
{\bf Lemma 8.1.}
{\it Let} $\hat v_i$ {\it satisfy} (7.1) {\it and} $\|\hat v_i\|_{\mu}<C$,
{\it where} $i=1,2$. {\it Let}
$$0<\tau\le\tau_1(\mu,\mu_0,C,\delta),\ \rho\ge\rho_1(\mu,\mu_0,C,\delta),
\eqno(8.1)$$
{\it where} $\tau_1,\rho_1$ {\it are the constants of Section} 4 {\it of}
[20] {\it and where} $\delta=1/2$, $2\le\mu_0<\mu$. {\it Then}
$$|||H_2-H_1|||_{\rho,\tau,\mu_0}\le
2(|||H^0_2-H^0_1|||_{\rho,\tau,\mu_0}+
|||Q^2_{\rho,\tau}-Q^1_{\rho,\tau}|||_{\rho,\tau,\mu_0}),\eqno(8.2)$$
{\it where} $H_i,\ H_i^0,\ Q^i_{\rho,\tau}$ {\it are the functions of}
(7.4), (7.7), (7.8), (7.14), (7.15) {\it for} $v=v_i$, $i=1,2$,
$|||\cdot|||_{\rho,\tau,\mu_0}$ {\it is defined as in} (7.5).
 {\it In addition,}
$$|||Q^2_{\rho,\tau}-Q^1_{\rho,\tau}|||_{\rho,\tau,\mu_0}\le
{24 c_5(\mu_0)C^2\over (1+2\tau\rho)^{\mu-\mu_0}}.\eqno(8.3)$$
In connection with (8.1) we remind that $\tau_1\in ]0,1[$ is sufficiently
small and $\rho_1$ is sufficiently great, see [20].

Lemma 8.1 follows from estimates mentioned as estimates (7.14), (7.15) of the
present paper (see estimates (4.3), (5.20), (5.22) of [20]) and from
Lemmas A.1, A.2, A.3 and estimate (A.5) (see Appendix to the present paper).

\vskip 2 mm
{\bf Lemma 8.2.}
{\it Let} $\hat v_i$ {\it satisfy} (7.1) {\it and} $\|\hat v_i\|_{\mu}<C$,
{\it where} $i=1,2$. {\it Let}
$$\eta(C,\rho_0,\mu)\le 1/2,\ \ \ln\rho_0\ge 2,\eqno(8.4)$$
{\it where} $\eta$ {\it is defined in} (7.6). {\it Let}
$$0<\tau_0<1,\ \ 2\le\mu_0<\mu,\ \ \rho=2\rho_0,\ \ \tau=\tau_0/2.$$
{\it Then}
$$\eqalignno{
&|||H^0_2-H^0_1|||_{\rho,\tau,\mu_0}\le (c_6+4c_7(\mu_0,\tau_0,\rho_0)C)
\Delta_{\rho_0,\rho,\mu_0},&(8.5)\cr
&\Delta_{\rho_0,\rho,\mu_0}\buildrel \rm def \over =
|||\chi_{\rho_0,\rho}(H_2-H_1)|||_{\rho_0,\mu_0},&(8.6)\cr}$$
{\it where} $H_i,\ H_i^0$ {\it are the functions of}
(4.13), (7.8) {\it for} $v=v_i$, $i=1,2$, $\chi_{\rho_0,\rho}$
{\it is the characteristic function of}
$\bar\Omega_{\rho}\b\bar\Omega_{\rho_0}$, $|||\cdot|||_{\rho,\tau,\mu_0}$
{\it in} (8.5) {\it is defined as in} (7.5),
$|||\cdot|||_{\rho_0,\mu_0}$  {\it in} (8.6) {\it is defined as in} (A.12) of
Appendix, $c_6$  {\it is defined by} (9.9), $c_7$ {\it is the constant} $c_8$
 {\it of} [20] ({\it that is}
$c_7(\mu,\tau,\rho)=3b_1(\mu)\tau^2+4b_2(\mu)\rho^{-1}+4b_3(\mu)\tau$,
{\it where} $b_1$, $b_2$, $b_3$ {\it are the constants of} [20]).

Lemma 8.2 is proved in Section 9.

\vskip 2 mm
{\bf Lemma 8.3.}
{\it Let the assumptions of Theorem} 2.1 {\it hold (for} $d=3$). {\it Let}
$$\rho_0\ge r_1(N,D,m,\sigma)\ \ {\it for\ some}\ \ \sigma>1,\eqno(8.7)$$
{\it where} $r_1$ {\it is the number of} (4.9). {\it Let}
$0<\tau_0<1$, $0<\mu_0$, $\rho=2\rho_0$, $\tau=\tau_0/2$. {\it Then}
$$\Delta_{\rho_0,\rho,\mu_0}\le
c_8\sigma^2e^{2\rho L}\|\Phi_2-\Phi_1\|(1+2\rho)^{\mu_0},\eqno(8.8)$$
{\it where} $\Delta_{\rho_0,\rho,\mu_0}$ {\it is defined by} (8.6),
$$c_8=(2\pi)^{-d}\int\limits_{\pa D}dx,\ \ L=\max\limits_{x\in\pa D}|x|,
\eqno(8.9)$$
$\|\Phi_2-\Phi_1\|$ {\it is defined according to} (2.4).

Lemma 8.3 is proved in Section 9.

\vskip 2 mm
{\bf Lemma 8.4.}
{\it Let the assumptions of Theorem} 2.1 {\it hold (for} $d=3$). {\it Let}
$$0<\tau\le\tau_2(m,\mu_0,N),\ \ \rho\ge\rho_2(m,\mu_0,N,D,\sigma),
\eqno(8.10)$$
{\it where} $2\le\mu_0<m$, $\sigma>1$, {\it and} $\tau_2$, $\rho_2$
{\it are constants such that} (8.10) {\it implies that}
$$\eqalignno{
&\tau\le\tau_1(m,\mu_0,c_4(m,3)N,1/2),\ \
\rho\ge\rho_1(m,\mu_0,c_4(m,3)N,1/2),&(8.11a)\cr
&\tau<1/2,\ \eta(c_4(m,3)N,\rho/2,m)\le 1/2,\ \ln(\rho/2)\ge 2,&(8.11b)\cr
&\rho/2\ge r_1(N,D,m,\sigma),&(8.11c)\cr}$$
{\it where}
$\tau_1$, $\rho_1$, $\eta$, $r_1$ {\it are the same that in}
(8.1), (8.4), (8.7), $c_4$ {\it is the constant of} (7.3). {\it Then}
$$|||H_2-H_1|||_{\rho,\tau,\mu_0}\le
c_9(N,D,m,\mu_0,\sigma,\tau)e^{2L\rho}\rho^{\mu_0}\|\Phi_2-\Phi_1\|+
c_{10}(N,m,\mu_0,\tau)\rho^{-(m-\mu_0)},\eqno(8.12)$$
{\it where}  $c_9$, $c_{10}$ {\it are some constants which can be given
explicitly}.

Lemma 8.4 follows from  formula (7.3) and Lemmas 8.1, 8.2, 8.3.

The final part of the proof of Theorem 2.2 for $d=3$ consists of the
following.
Using the inverse Fourier transform formula
$$v(x)=\int\limits_{\R^3}e^{-ipx}\hat v(p)dp,\ \ x\in\R^3,\eqno(8.13)$$
we have that
$$\eqalign{
&\|v_1-v_2\|_{L^{\infty}(D)}\le\sup\limits_{x\in\bar D}
\big|\int\limits_{\R^3}e^{-ipx}(\hat v_2(p)-\hat v_1(p))dp\big|\le\cr
&I_1(r)+I_2(r)\ \ {\rm for\ any}\ \ r>0,\cr}\eqno(8.14)$$
where
$$\eqalign{
&I_1(r)=\int\limits_{|p|\le r}|\hat v_2(p)-\hat v_1(p)|dp,\cr
&I_2(r)=\int\limits_{|p|\ge r}|\hat v_2(p)-\hat v_1(p)|dp.\cr}\eqno(8.15)$$

Using (7.3), (7.5), (7.14) for $H$, (7.16) we obtain that
$$|\hat v_2(p)-\hat v_1(p)|\le |||H_2-H_1|||_{\rho,\tau,2}(1+|p|)^{-2}
\eqno(8.16)$$
for $|p|<2\tau\rho$, at least, under the assumptions of Lemma 8.4 for
$\mu_0=2$.

Using (7.3) we obtain that
$$|\hat v_2(p)-\hat v_1(p)|\le 2c_4(m,3)N(1+|p|)^{-m},\ \ p\in\R^3.\eqno(8.17)
$$
Using (8.14)-(8.17) we have that
$$\eqalign{
&\|v_1-v_2\|_{L^{\infty}(D)}\le I_1(2\tau\rho)+I_2(2\tau\rho)\le\cr
&|||H_2-H_1|||_{\rho,\tau,2}\int\limits_{|p|\le 2\tau\rho}
{dp\over (1+|p|)^2}+
2c_4(m,3)N\int\limits_{|p|\ge 2\tau\rho}
{dp\over (1+|p|)^m}\le\cr
&8\pi\tau\rho |||H_2-H_1|||_{\rho,\tau,2}+
{8\pi c_4(m,3)N\over (m-3)(2\tau)^{m-3}}{1\over \rho^{m-3}},\cr}\eqno(8.18)$$
at least, under the assumptions of Lemma 8.4 for $\mu_0=2$.

Under the assumptions of Lemma 8.4 for $\mu_0=2$, using (8.18) and (8.12) we
obtain that
$$\|v_1-v_2\|_{L^{\infty}(D)}\le c_{11}(N,D,m,\sigma,\tau)e^{2L\rho}\rho^3
\|\Phi_2-\Phi_1\|+c_{12}(N,m,\tau)\rho^{-(m-3)},\eqno(8.19)$$
where $c_{11}$, $c_{12}$ are related in a simple way with $c_9$, $c_{10}$
for $\mu_0=2$.

Let now
$$\alpha\in ]0,1[,\ \beta={{1-\alpha}\over 2L},\ \delta=\|\Phi_1-\Phi_2\|,\
\rho=\beta\ln(3+\delta^{-1}),\eqno(8.20)$$
where $\delta$ is so small that $\rho\ge\rho_2(m,2,N,D,\sigma)$, where
$\rho_2$ is the constant of (8.10). Then, due to (8.19),
$$\eqalign{
&\|v_1-v_2\|_{L^{\infty}(D)}\le c_{11}(N,D,m,\sigma,\tau)
(3+\delta^{-1})^{2L\beta}(\beta\ln(3+\delta^{-1}))^3\delta+\cr
&c_{12}(N,D,m,\tau)(\beta\ln(3+\delta^{-1}))^{-(m-3)}=\cr
&c_{11}(N,D,m,\sigma,\tau)\beta^3(1+3\delta)^{1-\alpha}\delta^{\alpha}
(\ln(3+\delta^{-1}))^3+\cr
&c_{12}(N,D,m,\tau)\beta^{-(m-3)}
(\ln(3+\delta^{-1}))^{-(m-3)},\cr}\eqno(8.21)$$
where $\sigma$, $\tau$ are the same that in (8.10) for $\mu_0=2$ and
where $\alpha$, $\beta$ and $\delta$ are the same that in (8.20).

Using (8.21) we obtain that
$$\|v_1-v_2\|_{L^{\infty}(D)}\le c_{13}(N,D,m)
(\ln(3+\|\Phi_1-\Phi_2\|^{-1}))^{-(m-3)}\eqno(8.22)$$
for $\delta=\|\Phi_1-\Phi_2\|\le\delta_0(N,D,m)$, where $\delta_0$ is
sufficiently small positive constant. Estimate (8.22) in the general case
(with modified $c_{13}$) follows from (8.22) for

\noindent
$\delta=\|\Phi_1-\Phi_2\|\le\delta_0(N,D,m)$  and the property that
$\|v_j\|_{L^{\infty}(D)}\le c_{14}(m)N$ (for $d=3$).

Thus, Theorem 2.2 for $d=3$ is proved.

\vskip 2 mm
{\bf 9. Proofs of Lemmas 8.2 and 8.3}
\vskip 2 mm
{\it Proof of Lemma 8.2.}
Using the maximum principle for holomorphic functions it is sufficient to
prove that
$$\eqalign{
&\sup\limits_{(\lambda,p)\in b\Lambda_{\rho,\tau,\nu}^{\pm}}
(1+|p|)^{\mu_0}|H_2^0(\lambda(1\mp 0),p)-H_1^0(\lambda(1\mp 0),p)|\le\cr
&(c_6+4c_7(\mu_0,\tau_0,\rho_0)C)\Delta_{\rho_0,\rho,\mu_0},\cr}\eqno(9.1)$$
where $b\Lambda_{\rho,\tau,\nu}^+$, $b\Lambda_{\rho,\tau,\nu}^-$ are defined
in (6.17) (and where $H_i^0(\lambda(1-0),p)$, $i=1,2$, are considered for
$(\lambda,p)\in b\Lambda_{\rho,\tau,\nu}^+$,
$H_i^0(\lambda(1+0),p)$, $i=1,2$, are considered for
$(\lambda,p)\in b\Lambda_{\rho,\tau,\nu}^-$).

Using (7.8) and the Sohotsky-Plemelj formula, we have that
$$\eqalignno{
&H^0(\lambda(1-0),p)={1\over 2\pi i}\int\limits_{{\cal T}^+_{\rho/|p|}}
H(\zeta,p){d\zeta\over {\zeta-\lambda(1+0)}}+H(\lambda,p),\
(\lambda,p)\in b\Lambda_{\rho,\tau,\nu}^+,&(9.2a)\cr
&H^0(\lambda(1+0),p)=-{1\over 2\pi i}\int\limits_{{\cal T}^-_{\rho/|p|}}
H(\zeta,p){\lambda d\zeta\over \zeta(\zeta-\lambda(1-0))}+H(\lambda,p),\
(\lambda,p)\in b\Lambda_{\rho,\tau,\nu}^-.&(9.2b)\cr}$$
In addition, using the Cauchy-Green formula we have that
$$\eqalignno{
&H(\lambda,p)=-{1\over 2\pi i}\int\limits_{{\cal T}^+_{\rho/|p|}}
H(\zeta,p){d\zeta\over {\zeta-\lambda(1+0)}}+
{1\over 2\pi i}\int\limits_{{\cal T}^+_{\rho/|p|}}
H(\zeta,p){d\zeta\over {\zeta-\lambda}}-&(9.3a)\cr
&{1\over \pi}
\int\!\!\!\!\!\!\int\limits_{{\cal D}^+_{\rho_0/|p|}\b {\cal D}^+_{\rho/|p|}}
{\pa H(\zeta,p)\over \pa\bar\zeta}{dRe\zeta dIm\zeta\over {\zeta-\lambda}},\
(\lambda,p)\in b\Lambda_{\rho,\tau,\nu}^+,\cr
&H(\lambda,p)={1\over 2\pi i}\int\limits_{{\cal T}^-_{\rho/|p|}}
H(\zeta,p){\lambda d\zeta\over \zeta(\zeta-\lambda(1-0))}-
{1\over 2\pi i}\int\limits_{{\cal T}^-_{\rho/|p|}}
H(\zeta,p){\lambda d\zeta\over \zeta(\zeta-\lambda)}-&(9.3b)\cr
&{1\over \pi}
\int\!\!\!\!\!\!\int\limits_{{\cal D}^-_{\rho_0/|p|}\b {\cal D}^-_{\rho/|p|}}
{\pa H(\zeta,p)\over \pa\bar\zeta}
{\lambda dRe\zeta dIm\zeta\over \zeta(\zeta-\lambda)},\
(\lambda,p)\in b\Lambda_{\rho,\tau,\nu}^-\cr}$$
(where  the integrals along ${\cal T}^{\pm}_r$ are taken in the counter-
clockwise direction). In addition (see equation (A.11) of Appendix),
$${\pa H(\zeta,p)\over \pa\bar\zeta}=\{H,H\}(\zeta,p),\
(\zeta,p)\in\Lambda_{\rho_0,\tau_0,\nu},\eqno(9.4)$$
where $H(\zeta,p)$ in the left-hand side of (9.4) is defined according to
(7.4), $H$ in the right-hand side of (9.4) is the function of (4.13),
$\{\cdot,\cdot\}$ is  defined by (7.11).

Using (9.2), (9.3) we obtain that
$$\eqalignno{
&H^0(\lambda(1-0),p)={1\over 2\pi i}\int\limits_{{\cal T}^+_{\rho_0/|p|}}
H(\zeta,p){d\zeta\over {\zeta-\lambda}}-&(9.5a)\cr
&{1\over \pi}
\int\!\!\!\!\!\!\int\limits_{{\cal D}^+_{\rho_0/|p|}\b {\cal D}^+_{\rho/|p|}}
{\pa H(\zeta,p)\over \pa\bar\zeta}{dRe\zeta dIm\zeta\over {\zeta-\lambda}},\
(\lambda,p)\in b\Lambda_{\rho,\tau,\nu}^+,\cr
&H^0(\lambda(1+0),p)=-{1\over 2\pi i}\int\limits_{{\cal T}^-_{\rho_0/|p|}}
H(\zeta,p){\lambda d\zeta\over \zeta(\zeta-\lambda)}-&(9.5b)\cr
&{1\over \pi}
\int\!\!\!\!\!\!\int\limits_{{\cal D}^-_{\rho_0/|p|}\b {\cal D}^-_{\rho/|p|}}
{\pa H(\zeta,p)\over \pa\bar\zeta}
{\lambda dRe\zeta dIm\zeta\over \zeta(\zeta-\lambda)},\
(\lambda,p)\in b\Lambda_{\rho,\tau,\nu}^-.\cr}$$
Using (9.5), (9.4) for $H^0=H^0_i$, $H=H_i$, $i=1,2$, we obtain that:
$$\eqalignno{
&(H_2^0-H_1^0)(\lambda(1-0),p)=A^+(\lambda,p)+B^+(\lambda,p),&(9.6a)\cr
&A^+(\lambda,p)={1\over 2\pi i}\int\limits_{{\cal T}^+_{\rho_0/|p|}}
(H_2-H_1)(\zeta,p){d\zeta\over {\zeta-\lambda}},\cr
&B^+(\lambda,p)=-{1\over \pi}
\int\!\!\!\!\!\!\int\limits_{{\cal D}^+_{\rho_0/|p|}\b {\cal D}^+_{\rho/|p|}}
\bigl(\{H_2-H_1,H_2\}+\{H_1,H_2-H_1\}\bigr)
(\zeta,p){dRe\zeta dIm\zeta\over {\zeta-\lambda}},\cr}$$
where $(\lambda,p)\in b\Lambda_{\rho,\tau,\nu}^+$;
$$\eqalignno{
&(H_2^0-H_1^0)(\lambda(1+0),p)=A^-(\lambda,p)+B^-(\lambda,p),&(9.6b)\cr
&A^-(\lambda,p)=-{1\over 2\pi i}\int\limits_{{\cal T}^-_{\rho_0/|p|}}
(H_2-H_1)(\zeta,p){\lambda d\zeta\over \zeta(\zeta-\lambda)},\cr
&B^-(\lambda,p)=-{1\over \pi}
\int\!\!\!\!\!\!\int\limits_{{\cal D}^-_{\rho_0/|p|}\b {\cal D}^-_{\rho/|p|}}
\bigl(\{H_2-H_1,H_2\}+\{H_1,H_2-H_1\}\bigr)(\zeta,p)
{\lambda dRe\zeta dIm\zeta\over \zeta(\zeta-\lambda)},\cr}$$
where $(\lambda,p)\in b\Lambda_{\rho,\tau,\nu}^-$.

Estimates (9.1) follow from formulas (9.6) and from the estimates
$$\eqalignno{
&|A^{\pm}(\lambda,p)|\le c_6(1+|p|)^{-\mu_0}\Delta_{\rho_0,\rho,\mu_0},\ \
(\lambda,p)\in b\Lambda_{\rho,\tau,\nu}^{\pm},&(9.7)\cr
&|B^{\pm}(\lambda,p)|\le 4c_7(\mu_0,\tau_0,\rho_0)C(1+|p|)^{-\mu_0}
\Delta_{\rho_0,\rho,\mu_0},\ \
(\lambda,p)\in b\Lambda_{\rho,\tau,\nu}^{\pm},&(9.8)\cr}$$
where
$$\eqalign{
&c_6=\sup\limits_{r\in ]1/2,+\infty[}{q(r)\over {q(r)-q(2r)}},\cr
&q(r)=2r\bigl(1-\bigl(1-{1\over 4r^2}\bigr)^{1/2}\bigr).\cr}\eqno(9.9)$$

Note that
$$0<c_6\le (2\sqrt{3}-3)^{-1},\eqno(9.10)$$
where $c_6$ is defined by (9.9). Estimate (9.10) follows from the formulas
$$\eqalignno{
&c_6={1\over {1-2\sigma}},\ \ \sigma=\sup\limits_{\tau\in ]0,1[}
{{1-(1-(1/4)\tau)^{1/2}}\over {1-(1-\tau)^{1/2}}},&(9.11)\cr
&(1-\tau)^{1/2}\le 1-(1/2)\tau,\ \ 1-(1/4)\tau\ge a(1-(1/4)\tau)+1-a,&(9.12)
\cr
&a=2(2-\sqrt{3}),\ \tau\in ]0,1[.\cr}$$

Estimates (9.7) follow from the property that
$$(\zeta,p)\in  b\Lambda_{\rho_0,\tau_0,\nu}\ \ {\rm if}\ \ \zeta\in
{\cal T}_{\rho_0/|p|}^{\pm},\ |p|<2\tau\rho,\ \rho=2\rho_0>0,\
\tau=\tau_0/2,\ 0<\tau_0<1,\eqno(9.13)$$
the properties that
$$\eqalign{
&H_i\in {\cal C}(\Lambda_{\rho_0,\tau_0,\nu}\cup b\Lambda_{\rho_0,\tau_0,\nu})
,\ i=1,2,\cr
&|(H_2-H_1)(\lambda,p)|\le (1+|p|)^{-\mu_0}\Delta_{\rho_0,\rho,\mu_0} ,\ \
(\lambda,p)\in b\Lambda_{\rho_0,\tau_0,\nu}\cr}\eqno(9.14)$$
(see formulas (A.8), (A.9), (A.12), (6.19), (7.4))
and from the formulas
$${1\over 2\pi}\int\limits_{{\cal T}_r^-}{|\lambda| |d\zeta|\over
|\zeta| |\zeta-\lambda|}=
{1\over 2\pi}\int\limits_{{\cal T}_r^+}{|d\zeta^{\prime}|\over
|\zeta^{\prime}-\lambda^{\prime}|}\le {q(r)\over {q(r)-q(2r)}},\eqno(9.15)$$
$\lambda\in {\cal T}_{2r}^-$, $\lambda^{\prime}\in {\cal T}_{2r}^+$,
$r>1/2$. In turn, formulas (9.15) follow from the property that
$z^{-1}\in {\cal T}_r^+$ if $z\in {\cal T}_r^-$ and from the formula that
$q(r)$ is the radius of ${\cal T}_r^+$, where $r>1/2$.

Estimates (9.8) follow from formulas (A.8), (A.9) (for $\rho=\rho_0$),
formula (8.4),

\noindent
Lemma A.5 (for $\rho=\rho_0$, $\mu=\mu_0$), Lemma A.6
(for $\rho=\rho_0$) and the property
$$\lambda\in {\cal T}_{\rho/|p|}^{\pm}\subset {\cal D}_{\rho_0/|p|}^{\pm}\ \
{\rm if}\ \ (\lambda,p)\in b\Lambda_{\rho,\tau,\nu}^{\pm},\ \rho=2\rho_0>0,\
\tau=\tau_0/2,\ 0<\tau_0<1.\eqno(9.16)$$

Lemma 8.2 is proved.

\vskip 2 mm
{\it Proof of Lemma 8.3.}
Using (A.12), the formulas
$$\Omega\approx\Theta,\ \ \Omega_{\rho}\approx\Theta_{\rho},\
\bar\Omega_{\rho}\approx\bar\Theta_{\rho}$$
(see  (6.4)), and formulas (6.2), (4.13),  we have that
$$\Delta_{\rho_0,\rho,\mu_0}\le\sup\limits_{\scriptstyle
(k,l)\in\bar\Theta_{\rho}\b\Theta_{\rho_0}}
(1+|k-l|)^{\mu_0}|h_2(k,l)-h_1(k,l)|,\eqno(9.17)$$
where $\Delta_{\rho_0,\rho,\mu_0}$ is defined by (8.6), $h_1$, $h_2$ are the
functions of Section 5.

Estimate (8.8) follows from formulas (9.17), (5.1), (4.7), (4.9), the property
 that $|k-l|\le 2\rho$  for $(k,l)\in\bar\Theta_{\rho}$
and the property that $|e^{ikx}|\le e^{\rho L}$,
$|e^{ilx}|\le e^{\rho L}$ for $k,l\in\bar\Theta_{\rho}$, $x\in\pa D$.

Lemma 8.3 is proved.

\vskip 2 mm
{\bf Appendix}

In this appendix we recall some results of [20] used in Sections 8 and 9 of
the present work.

Consider the equation
$$U=U^0+M_{\rho,\tau}(U),\ \ \rho>0,\ \tau\in ]0,1[,\eqno(A.1)$$
for unknown $U$, where $U^0$, $U$ are functions on $\Lambda_{\rho,\tau,\nu}$
defined in (6.13), $M_{\rho,\tau}$ is the map defined by (7.9).

In particular, (7.7) can be written as (A.1) with $U=H$,
$U^0=H^0+Q_{\rho,\tau}$.

\vskip 2 mm
{\bf Lemma A.1.}
{\sl
Let $\rho>0$, $\nu\in\S^2$, $\tau\in ]0,1[$, $\mu\ge 2$ and
$0<r<(2c_7(\mu,\tau,\rho))^{-1}$, where $c_7$ is the constant $c_8$ of [20]
(that is $c_7(\mu,\tau,\rho)=3b_1(\mu)\tau^2+4b_2(\mu)\rho^{-1}+4b_3(\mu)\tau$
, where $b_1$, $b_2$, $b_3$ are the constants of [20]). Let
$U^0\in L^{\infty}(\Lambda_{\rho,\tau,\nu})$ and
$|||U^0|||_{\rho,\tau,\mu}\le r/2$ (see definition (7.5)).
Then equation (A.1) is uniquely solvable for
$U\in L^{\infty}(\Lambda_{\rho,\tau,\nu})$,
$|||U|||_{\rho,\tau,\mu}\le r$, and $U$ can be found by the method of
successive approximations, in addition,
$$|||U-(M_{\rho,\tau,U^0})^n(0)|||_{\rho,\tau,\mu}\le
{r(2c_7(\mu,\tau,\rho)r)^n\over 2(1-2c_7(\mu,\tau,\rho)r)},\ n\in\N,\eqno(A.2)
$$
where $M_{\rho,\tau,U^0}$ denotes the map $U\to U^0+M_{\rho,\tau}(U)$.
}

Lemma A.1 corresponds to Lemma 4.4 of [20].

\vskip 2 mm
{\bf Lemma A.2.}
{\sl
Let the assumptions of Lemma A.1 be fulfilled. Let also

\noindent
$\tilde U^0\in L^{\infty}(\Lambda_{\rho,\tau,\nu})$,
$|||\tilde U^0|||_{\rho,\tau,\mu}\le r/2$, and $\tilde U$ denote the solution
of (A.1) with $U^0$  replaced by $\tilde U^0$, where
$\tilde U\in L^{\infty}(\Lambda_{\rho,\tau,\nu})$,
$|||\tilde U|||_{\rho,\tau,\mu}\le r$. Then
$$|||U-\tilde U|||_{\rho,\tau,\mu}\le (1-2c_7(\mu,\tau,\rho)r)^{-1}
|||U^0-\tilde U^0|||_{\rho,\tau,\mu}.\eqno(A.3)$$
}

Lemma A.2 corresponds to Lemma 4.5 of [20].

Let
$$r_{min}\buildrel \rm def \over ={2C\over {1-\eta(C,\rho,\mu)}}+
{2b_4(\mu_0)C^2\over (1-\eta(C,\rho,\mu))^2}
\biggl(1+{3\over (1+2\tau\rho)^{\mu-\mu_0}}\biggr),\eqno(A.4)$$
where $C>0$, $\rho\ge e^2$, $2\le\mu_0<\mu$, $\tau\in ]0,1[$,
$\eta$ is defined as in (7.6), $b_4$ is the constant of [20] (that is
$b_4(\mu)=\pi^{-1}(b_1(\mu)n_1+b_2(\mu)n_2+b_3(\mu)n_3)$,
where $b_1$, $b_2$, $b_3$ are the constants of [20] and $n_1$, $n_2$,
$n_3$ are the constants of Lemma 11.1 of [19]).
In addition, from estimates (7.14), (7.15) of the present paper (see
estimates (3.3), (4.20), (4.22) of [20]) it follows that
$$\eqalign{
&|||H|||_{\rho,\tau,\mu_0}\le r_{min}\cr
&|||H^0|||_{\rho,\tau,\mu_0}+|||Q_{\rho,\tau}|||_{\rho,\tau,\mu_0}
\le r_{min}/2.\cr}\eqno(A.5)$$

\vskip 2 mm
{\bf Lemma A.3.}
{\sl
Let $C>0$, $2\le\mu_0<\mu$, $0<\delta<1$,
$$0<\tau\le\tau_1(\mu,\mu_0,C,\delta),\ \ \rho\ge\rho_1(\mu,\mu_0,C,\delta),
\eqno(A.6)$$
where $\tau_1$, $\rho_1$ are the constants of (4.36) of [20] (in particular,
$\tau_1\in ]0,1[$ is sufficiently small and $\rho_1$ is sufficiently great).
Then
$$\eqalign{
&0\le\eta(C,\rho,\mu)<\delta,\cr
&0\le 2c_7(\mu_0,\tau,\rho)r_{min}(\mu,\mu_0,\tau,\rho,C)<\delta.\cr}
\eqno(A.7)$$
}

Lemma A.3 corresponds to estimates (4.36) of [20].

\vskip 2 mm
{\bf Lemma A.4.}
{\sl
Let $v$ satisfy (7.1), $\|\hat v\|_{\mu}\le C$, $\rho$ satisfy (7.6). Let $H$
be the function of (4.13). Then:
$$\eqalignno{
&H\in {\cal C}(\Omega\b\Omega_{\rho}),&(A.8)\cr
&|H(k,p)|\le {C\over (1-\eta(C,\rho,\mu))(1+|p|)^{\mu}},\ \
(k,p)\in\Omega\b\Omega_{\rho},&(A.9)\cr
&\hat v(p)=\lim_{\scriptstyle |k|\to\infty,\atop (k,p)\in\Omega}H(k,p),\
p\in\R^3,&(A.10)\cr}$$
where ${\cal C}$ denotes the space of continuous functions, $\Omega$,
$\Omega_{\rho}$ are defined in (4.14), (6.3), $\eta$ is defined as in (7.6),
$|k|=((Re\,k)^2+(Im\,k)^2)^{1/2}$;
$${\pa\over \pa\bar\lambda}H(k(\lambda,p),p)=\{H,H\}(\lambda,p),\ \
(\lambda,p)\in\Lambda_{\rho,\nu},\eqno(A.11)$$
where $k(\lambda,p)$ is defined in (6.11b), $\Lambda_{\rho,\nu}$ is defined
in (6.13), $\{\cdot,\cdot\}$ is defined by (7.11).
}

Lemma A.4 corresponds to formulas (3.2)-(3.4), (3.22) of [20]. Let
$$\eqalign{
&L^{\infty}_{\mu}(\Omega\b\bar\Omega_{\rho})=
\{U\in L^{\infty}(\Omega\b\bar\Omega_{\rho}):\ |||U|||_{\rho,\mu}< +\infty\},
\cr
&|||U|||_{\rho,\mu}=ess\ \sup\limits_{(k,p)\in\Omega\b\bar\Omega_{\rho}}
(1+|p|)^{\mu}|U(k,p)|,\ \rho>0,\ \mu>0,\cr}\eqno(A.12)$$
where $\Omega$, $\bar\Omega_{\rho}$ are defined in (4.14), (6.3).

\vskip 2 mm
{\bf Lemma A.5.}
{\sl
Let $U_1,U_2\in L^{\infty}_{\mu}(\Omega\b\bar\Omega_{\rho})$, where
$\rho>0$, $\mu\ge 2$. Let $\{\cdot,\cdot\}$ be defined by (7.11). Then
$$\{U_1,U_2\}\in L^{\infty}_{local}(\Lambda_{\rho,\nu}),\eqno(A.13)$$
$$\eqalign{
&|\{U_1,U_2\}(\lambda,p)|\le {|||U_1|||_{\rho,\mu}|||U_2|||_{\rho,\mu}\over
(1+|p|)^{\mu}}b(\mu,|\lambda|,|p|)\cr
&{for\ almost\ all}\ \ (\lambda,p)\in\Lambda_{\rho,\nu},\
b(\mu,|\lambda|,|p|)=\cr
&\biggl({b_1(\mu)|\lambda|\over (|\lambda|^2+1)^2}+
{b_2(\mu)|p|||\lambda|^2-1|\over
 |\lambda|^2(1+|p|(|\lambda|+|\lambda|^{-1}))^2}+
{b_3(\mu)|p|\over
 |\lambda|(1+|p|(|\lambda|+|\lambda|^{-1}))}\biggl),\cr}\eqno(A.14)$$
where $\Lambda_{\rho,\nu}$ is defined in (6.13), $b_1$, $b_2$, $b_3$ are
the constants of [20].
}

Lemma A.5 corresponds to estimates (3.26), (3.27) of [20].

Consider
$$\eqalign{
&u_1(\zeta,s)={|\zeta|\over (|\zeta|^2+1)^2},\ \
u_2(\zeta,s)={(|\zeta|^2+1)s\over |\zeta|^2(1+s(|\zeta|+|\zeta|^{-1}))^2},\cr
&u_3(\zeta,s)={s\over |\zeta|(1+s(|\zeta|+|\zeta|^{-1}))},\cr}\eqno(A.15)$$
where $\zeta\in\C$, $s>0$.

Note that
$$0<b(\mu,|\lambda|,|p|)\le\sum_{j=1}^3b_j(\mu)u_j(\lambda,|p|),\eqno(A.16)$$
where $b$ is the function of (A.14), $(\lambda,p)\in\Lambda_{\rho,\nu}$,
$\rho>0$, $\mu\ge 2$.

\vskip 2 mm
{\bf Lemma A.6.}
{\sl
Let $u_1,u_2,u_3$ be defined by (A.15). Then:
$$\int\limits_{{\cal D}_r^-}u_j(\zeta,s){|\lambda|\over |\zeta|}
{d Re\,\zeta d Im\,\zeta\over |\zeta-\lambda|}=
\int\limits_{{\cal D}_r^+}u_j(\zeta,s)
{d Re\,\zeta d Im\,\zeta\over |\zeta-\lambda^{-1}|},\eqno(A.17)$$
$j=1,2,3$, $r>1/2$, $\lambda\in {\cal D}_r^-$, $s>0$,
$$\int\limits_{{\cal D}^+_{\rho/|p|}}u_j(\zeta,|p|)
{d Re\,\zeta d Im\,\zeta\over |\zeta-\lambda|}\le\left\{\matrix{
(3/4)\pi (|p|/\rho)^2\ &\ {for}\ j=1,\hfill\cr
4\pi/\rho\ &\ {for}\ j=2,\hfill\cr
2\pi |p|/\rho\ &\ {for}\ j=3,\hfill\cr}\right.\eqno(A.18)$$
$0<|p|<2\tau\rho$, $0<\tau<1$, $\lambda\in {\cal D}^+_{\rho/|p|}$,
where ${\cal D}_r^{\pm}$ are defined in (6.18).
}

Lemma A.6 corresponds to formulas (7.2), (7.10) of [20].

\vskip 4 mm
{\bf References}
\vskip 2 mm
\item{[  1]} G.Alessandrini, {\it Stable determination of conductivity
by boundary measurements}, Appl. Anal. {\bf 27} (1988), 153-172.
\item{[  2]} R.Beals and R.R.Coifman, {\it Multidimensional inverse
scattering and nonlinear partial differential equations}, Proc. Symp. Pure
Math. {\bf 43} (1985), 45-70.
\item{[   3]} R.M.Brown and G.Uhlmann, {\it Uniqueness in the inverse
conductivity problem for nonsmooth conductivities in two dimensions},
Comm. Partial Diff. Eq. {\bf 22} (1997), 1009-1027.
\item{[   4]} A.L.Bukhgeim, {\it Recovering a potential from Cauchy
data in the two-dimensional case}, J.Inverse Ill-Posed Probl. {\bf 16}(1)
(2008), 19-33.
\item{[   5]} A.-P.Calder\'on, {\it On an inverse boundary value problem},
Seminar on Numerical Analysis and its Applications to Continuum Physics
(Rio de Janeiro, 1980), pp.65-73, Soc. Brasil. Mat. Rio de Janeiro, 1980.
\item{[   6]} L.D.Druskin, {\it The unique solution of the inverse problem
in electrical surveying and electrical well logging for piecewise-constant
conductivity}, Physics of the Solid Earth {\bf 18}(1) (1982), 51-53.
\item{[   7]} L.D.Faddeev, {\it Growing solutions of the Schr\"odinger
equation}, Dokl. Akad. Nauk SSSR {\bf 165} (1965), 514-517 (in Russian);
English Transl.: Sov. Phys. Dokl. {\bf 10} (1966), 1033-1035.
\item{[   8]} L.D.Faddeev, {\it Inverse problem of quantum scattering theory
II}, Itogi Nauki i Tekhniki, Sovr. Prob. Math. {\bf 3} (1974), 93-180
(in Russian); English Transl.: J.Sov. Math. {\bf 5} (1976), 334-396.
\item{[   9]} I.M.Gelfand, {\it Some problems of functional analysis and
algebra}, Proceedings of the International Congress of Mathematicians,
Amsterdam, 1954, pp.253-276.
\item{[  10]} G.M.Henkin and R.G.Novikov, {\it The $\bar\pa$- equation in the
multidimensional inverse scattering problem}, Uspekhi Mat. Nauk {\bf 42(3)}
(1987), 93-152 (in Russian); English Transl.: Russ. Math. Surv. {\bf 42(3)}
(1987), 109-180.
\item{[  11]} R.Kohn and M.Vogelius, {\it Determining conductivity by
boundary measurements II, Interior results}, Comm. Pure Appl. Math. {\bf 38}
(1985), 643-667.
\item{[  12]} M.Lassas, J.L.Mueller and S.Siltanen, {\it Mapping properties
of the nonlinear Fourier transform in dimension two}, Comm.Partial
Differential Equations {\bf 32}(4-6) (2007), 591-610.
\item{[  13]} N.Mandache, {\it Exponential instability in an inverse problem
for the Schr\"odinger equation}, Inverse Problems {\bf 17} (2001), 1435-1444.
\item{[  14]} A.I.Nachman, {\it Reconstructions from boundary measurements},
Ann. Math. {\bf 128} (1988), 531-576.
\item{[  15]} A.I.Nachman, {\it Global uniqueness for a two-dimensional
inverse boundary value problem}, Ann, Math. {\bf 142} (1995), 71-96.
\item{[  16]} R.G.Novikov, {\it Multidimensional inverse spectral problem
for the equation $-\Delta\psi+(v(x)-Eu(x))\psi=0$}, Funkt. Anal. i Pril.
{\bf 22(4)} (1988), 11-22 (in Russian); English Transl.: Funct. Anal. and
Appl. {\bf 22} (1988), 263-272.
\item{[  17]} R.G.Novikov, {\it Formulae and equations for finding scattering
data from the Dirichlet-to-Neumann map with nonzero background potential},
Inverse Problems {\bf 21} (2005), 257-270.
\item{[  18]} R.G.Novikov, {\it The $\bar\pa$- approach to approximate inverse
scattering at fixed energy in three dimensions}, International Mathematics
Research Papers, {\bf 2005:6}, (2005), 287-349.
\item{[  19]} R.G.Novikov, {\it On non-overdetermined inverse scattering
at zero energy in three dimensions}, Ann. Scuola Norm. Sup. Pisa Cl. Sci.
{\bf 5} (2006), 279-328
\item{[  20]} R.G.Novikov, {\it An effectivization of the global
reconstruction in the Gel'fand-Calderon inverse problem in three dimensions},
Contemporary Mathematics, {\bf 494} (2009), 161-184.
\item{[  21]} R.G.Novikov and N.N.Novikova, {\it On stable determination of
potential by boundary measurements}, ESAIM: Proceedings {\bf 26} (2009),
94-99.
\item{[  22]} V.P.Palamodov, {\it Gabor analysis of the continuum model
for impedance tomography}, Ark.Mat. {\bf 40} (2002), 169-187.
\item{[  23]} J.Sylvester and G.Uhlmann, {\it A global uniqueness theorem
for an inverse boundary value problem}, Ann. Math. {\bf 125} (1987),
153-169.
\item{[ 24]} R.Novikov and M.Santacesaria, {\it A global stability estimate
for the Gel'fand-Calderon inverse problem in two dimensions}, in preparation.

\end